\documentclass[11pt,tgrind,psbox]{article}
\usepackage{latexsym}
\usepackage{graphics}
\newcommand {\be}[1]{\begin{equation}\label{#1}}
\newcommand {\ee}{\end{equation}}
\newcommand {\bea}{\begin{eqnarray}}
\newcommand {\eea}{\end{eqnarray}}

\newcommand{\proof}{{\bf Proof : }}
\newcommand{\qed}{\hfill $\Box$}
\newcommand{\remark}{{\bf Remark : }}
\newcommand{\remarks}{{\bf Remarks : }}
\newcommand{\vvert}{\vspace{.1in}}
\newcommand{\pr}{{\rm Prob}}
\newcommand{\cZ}{{\cal Z}_+^d}
\newcommand{\cZm}{{\cal Z}_+^{m+1}}
\newcommand{\cZmm}{{\cal Z}_+^{m+3}}
\newcommand{\cZLambda}{{\cal Z}_{\Lambda}}
\newtheorem{theorem}{Theorem}

\newtheorem{prop}{Proposition}

\newtheorem{Defi}{Definition}

\title{\bf Computing stationary probability distributions and large deviation rates for
constrained random walks. The undecidability results.}

\author{David Gamarnik
\thanks{IBM T.J. Watson Research Center, Yorktown Heights, NY 10598, USA.
        Email address: gamarnik@watson.ibm.com }
}

\begin{document}

\maketitle

\centerline{{\bf Keywords:} positive recurrence, computability,
Lyapunov function.}

\begin{abstract}
Our model is a constrained homogeneous random walk in $\cZ$. The
convergence to stationarity for such a random walk can often be
checked by constructing a Lyapunov function. The same Lyapunov
function can also be used for computing approximately the
stationary distribution of this random walk, using methods
developed by Meyn and Tweedie in \cite{metwee_paper}. In this
paper we show that, for this type of random walks, computing the
stationary probability exactly is an undecidable problem: no
algorithm can exist to achieve this task. We then prove that
computing large deviation rates for this model is also an
undecidable problem. We extend these results to a certain type of
queueing systems. The implication of these results is that no
useful formulas for computing stationary probabilities and large
deviations rates can exist in these systems.
\end{abstract}

\section{Introduction}\label{introduction}
The main model considered in this paper is a constrained
homogeneous random walk in a $d$-dimensional nonnegative orthant
$\cZ$, where $\cZ$ is the space of $d$-dimensional vectors with
integral nonnegative components. Specifically, the transitions
with positive probabilities can occur only to neighboring states
and the transition probabilities depend only on the face that the
current state of the random walk belongs to, but not on the size
of the components of the state.

Ever since the appearance of the papers by Malyshev
\cite{malyshev}, \cite{MalyshevAnalyticMethods},
\cite{MalyshevAsymptotic} and Menshikov \cite{menshikov}, random
walks in $\cZ$ have assumed a prominent role in modelling and
analysis of queueing networks of certain type, for example
Markovian queueing networks. Specifically, the question of
positive recurrence or stability was analyzed. One of the main
techniques used for the stability analysis of these type of random
walks is a Lyapunov function technique also known as Foster's
criteria. A comprehensive study of constrained random walks in
$\cZ$ was conducted by Fayolle, Malyshev and Menshikov in
\cite{malyshev_menshikov_book} and many additional results
appeared after the book was published. Specifically, a very
interesting connection between constrained random walks and
general dynamical systems on compact manifolds was established by
Malyshev \cite{malyshev_dynamics}. Exact conditions for positive
recurrence for the case $d\leq 4$ were obtained by Ignatyuk and
Malyshev in \cite{ign_malyshev}. The large deviation principle for
special cases and modifications of random walks in $\cZ$ was
established by Ignatyuk, Malyshev and Scherbakov  in
\cite{IgnatyukMalyshevScherbakov}. This followed by efforts to
actually compute the large deviation rates, which turned out to be
a very complicated problem. See for example Kurkova and Suhov
\cite{suhov}, where large deviations limits are computed for a
random walk in ${\cal Z}_+^2$ arising from joint-the-shortest
queueing system. The analysis uses a fairly complicated
complex-analytic techniques. The goal of the current paper is to
explain the difficulty in obtaining such results for general
dimensions.

Analysis of random walks arising from special types of multiclass
queueing networks became a subject of particularly aggressive
research efforts during the previous decade. Many interesting and
deep results were established which connect stability of such
queueing networks with stability of corresponding fluid models,
obtained by Law of Large Numbers type of rescaling. This research
direction was initiated in pioneering works by Rybko and Stolyar
\cite{rs} and Dai \cite{dai}, where it was shown that stability of
a fluid model implies stability of the underlying queueing system.
The converse of this result is not true in general, see Dai,
Hasenbein and Vande Vate \cite{daivvhasscounter}, Bramson
\cite{bramsoncounter}, but is true under some stronger conditions,
Dai \cite{dai_converse}, Meyn \cite{meyn_converse}. Despite these
results, to the day no full characterization of stable queueing
networks is available. Stability was characterized only for
special types of queueing networks or special scheduling polices.
For example feedforward networks are known to be stable for all
work-conserving policies, Down and Meyn \cite{medo}, Dai
\cite{dai}. Stability of fluid networks with two processing
stations operating under arbitrary work-conserving policies is
fully characterized in Bertsimas, Gamarnik and Tsitsiklis
\cite{bgt} by means of a linear programming and in Dai and Vande
Vate \cite{dai_vv} by direct methods. The question of computing
stationary distributions comes naturally after the question of
stability. Several results are available again in the context of
multiclass queueing networks. Some of the results were obtained
using quadratic Lyapunov functions, Bertsimas, Paschalidis and
Tsitsiklis \cite{bpt}, Kumar and Kumar \cite{kumar_bounds}, Kumar
and Meyn \cite{kumar_meyn}, Kumar and Morrison
\cite{kumarmorrison}, using piece-wise linear Lyapunov functions
in Bertsimas, Gamarnik and Tsitsiklis \cite{bgt_perf}, and using
more direct methods, Bertsimas and Nino-Mora
\cite{BertsimasNinomoraI}, \cite{BertsimasNinomoraII} All of the
results obtain only bounds on the stationary probabilities.
Computing exactly the stationary probabilities seems beyond the
existing techniques.

\section{Our results}
It was established by the author in \cite{gamarnik_decidability}
that positive recurrence of a constrained homogeneous random walk
in $\cZ$ is an undecidable property. Meaning, no algorithm can
possibly exists which given the description of the random walks
(given the dimension and the transition matrix) will be able to
check whether the walk is positive recurrent. This result was also
established for queueing systems operating under the class of so
called {\it generalized priority} policies. This result explains
the difficulty in stability analysis by stating that these
problems are simply insolvable. It was conjectured in the same
paper that the stability of multiclass  queueing networks
operating under the class of much studied priority or
First-In-First-Out policies, is undecidable as well. The
conjecture remains unproven.

In the current paper we continue the decidability analysis of
constrained random walks by asking the following question: given a
constrained homogeneous random walk, can we compute its stationary
distribution, provided that the existence of a stationary
distribution can be checked, for example, by constructing a
Lyapunov function? To put this question into a proper computation
theoretic framework, we ask the following question. Given a
constrained homogeneous random walk, which possesses a unique
stationary distribution $\pi$, given a state $q\in\cZ$, for
example $q=0$, and given a rational value $r>0$, is it true that
the stationary probability $\pi(q)$ of this state satisfies
$\pi(q)\leq r$? In this paper we  prove that this problem is
undecidable, \emph{ even if} a Lyapunov function witnessing
positive recurrence is available. Thus, no algorithm can exist
which given a positive recurrent constrained homogeneous random
walk computes its stationary distribution. Specifically, the
stationary distribution cannot be written down in any constructive
way using some formulas. Contrast this with random walks
corresponding to product form type networks, for example Jackson
networks, for which a very simple formula is available. We then
prove that computing large deviations limits for the same model is
an undecidable problem as well. In particular, we show that given
a random walk in $\cZ$ with a unique stationary distribution
(witnessed, for example, by a Lyapunov function), and given a
vector $v\in\Re^d$, the problem of deciding whether
$\lim_{n\rightarrow\infty}\log(\pi(vn))/n$ is finite, is
undecidable. We extend these results to queueing systems operating
under a class of generalized priority policies. Finally, we
observe that, nevertheless, estimating stationary distribution of
a constrained random walk is a decidable problem, if one is
willing to tolerate a  two-sided error and a Lyapunov function
exists. Specifically, given such a random walk with the unique
stationary distribution $\pi$, given a Lyapunov function, given a
state $q\in\cZ$  and any value $\epsilon>0$, an interval
$(r,r+\epsilon)$ can be constructed which contains $\pi(q)$. This
result is an easy consequence of a powerful result obtained by
Meyn and Tweedie \cite{metwee_paper}, which obtains exponential
bounds on the mixing rate of Markov chains, using Lyapunov
function methods. We note that such approximation result cannot be
obtained for large deviations rates since, as we mentioned above,
even determining whether a given large deviation rate is finite is
an undecidable problem.

The remainder of the paper is organized as follows. In the
following section we describe our model -- constrained homogeneous
random walk in $\cZ$ and introduce Lyapunov functions. In Section
\ref{CounterMachine} we introduce a counter machine -- a
modification of a Turing machine which for us is the main tool for
establishing the undecidability results. In Section
\ref{StationaryUndecidable} we prove that computing a stationary
distribution of a positive recurrent random walk in $\cZ$ is an
undecidable problem. In Section \ref{LDUndecidable} we prove that
computing large deviations rates for positive recurrent random
walks in $\cZ$ is an undecidable problem as well. Extension of
these results to queueing systems is established in Section
\ref{queues}. In Section \ref{StationaryDecidable} we show how
stationary distribution can be computed with a two-sided error
using a Lyapunov function technique and Meyn and Tweedie results
from \cite{metwee_paper}. Conclusions and open problems are
discussed in Section \ref{conclusions}.

\section{Constrained homogeneous random walk in $\cZ$. Lyapunov
function and stationary distribution}\label{definitions} Let $\cZ$
denote the space of $d$-dimensional vectors with nonnegative
integer components. Our model is a random walk
$Q(t),t=0,1,2,\ldots\,\,$ which has $\cZ$ as  a state-space. For
each $\Lambda\subset \{1,2,\ldots,d\}$ let $\cZLambda$ denote the
corresponding face:
\[ \cZLambda=\{(z_1,z_2,\ldots,z_d)\in \cZ:
z_i>0 \,{\rm for}\,i\in \Lambda, z_i=0 \,{\rm for}\,i\notin
\Lambda\}.
\]
The transition probabilities  are face-homogeneous -- they depend
entirely on the face the random walk is currently on. In addition
the transition vectors have at most unit length in $\max$ norm. In
other words, for each $\Lambda\subset \{1,2,\ldots,d\}$ and each
$\Delta \in \{-1,0,1\}^d$ a certain value $p(\Lambda,\Delta)$ (the
transition probability) is defined. These values satisfy
\[
\sum_{\Delta\in\{-1,0,1\}^d}p(\Lambda,\Delta)=1
\]
for each $\Lambda$ and $p(\Lambda,\Delta)=0$ if for some $i\notin
\Lambda, \Delta_i=-1$. The latter condition is simply a
consistency condition which prevents transitions into states with
negative components. Given a current state $Q(t)\in \cZ$ of the
random walk, the next state is chosen to be $Q(t)+\Delta $ with
probability $p(\Lambda,\Delta)$, if the state $Q(t)$ belongs to
the face $\cZLambda$. We will also write $p(q,q')$ instead of
$p(\Lambda,\Delta)$ if the state $q\in\cZLambda$ and
$q'-q=\Delta$. We denote by $p^{(t)}(q,q')$ the $t$-step
transition probabilities: $\pr\{Q(t)=q'|Q(0)=q\}$. The model above
will be referred to as constrained homogeneous random walk in
$\cZ$. We will say that our walk is deterministic if
$p(\Lambda,\Delta)\in\{0,1\}$ for all $\Lambda$ and $\Delta$. In
other words, the transition vector $\Delta=\Delta(\Lambda)$
deterministically depends on the face. The set of parameters
$p(\Lambda,\Delta)$ is finite, and, in particular, it contains
$6^d$ elements corresponding to $2^d$ faces $\cZLambda$ and $3^d$
transition vectors $\Delta$ per face (with some transitions
occurring with zero probability). Let $||Q(t)||$ denote $L_1$
norm. That is $||Q(t)||=\sum_{i\leq d}Q_i(t)$.

For any state $q\in\cZ$ and subset $X\subset\cZ$, let $T=T(q,X)$
denote the first hitting time for the set $X$ when the initial
state of the walk is $q$, including the possibility $T=\infty$.
That is \be{eq:hittime}
 T=\min\{t:Q(t)\in X|Q(0)=q\}.
\ee
The following definition is standard in the theory or infinite
Markov chains.

\begin{Defi}\label{PositiveRecurrence}
A homogeneous random walk is defined to be positive recurrent or
stable if there exist some $C>0$  such that the random walk visits
the set $X_C\equiv\{z\in \cZ:\sum_{i=1}^d z_i\leq C\}$ infinitely
often with probability one, and $E[T(q,X_C)]$ is finite for all
$q\in\cZ$.
\end{Defi}

Stability of a constrained homogeneous random walk $Q(t)$ can be
checked, for example, by constructing a suitable Lyapunov
function.
\begin{Defi}\label{Lyapunov}  A function $\Phi:\cZ\rightarrow \Re_+$ is defined
to be a Lyapunov function with drift $-\gamma<0$ and exception set
${\cal B}\subset \cZ$ if $|{\cal B}|<\infty$ and  for every state
$q\notin {\cal B}$ \be{drift}
E[\Phi(Q(t+1)|Q(t)=q]-\Phi(q)=\sum_{q\in\cZ}\Phi(q')p(q,q')-\Phi(q)\leq-\gamma.
\ee
\end{Defi}

In other words, the expected value of the Lyapunov function should
decrease at each time step, whenever the random walk is outside of
the exception set. Existence of a Lyapunov function under some
additional assumptions implies stability. For a comprehensive
survey of Lyapunov function methods see Meyn and Tweedie
\cite{metwee}. Various forms of Lyapunov functions, specifically
piece-wise linear and quadratic Lyapunov functions were used to
prove stability of random walks corresponding to Markovian
queueing networks \cite{bgt_perf}, \cite{bpt}, \cite{kumarou},
\cite{kumar_meyn}, \cite{kumarmorrison}, \cite{kumar_bounds}. In
some simple cases even linear Lyapunov function of the form
$\Phi(q)=\sum_{i=1}^dw_iq_i$ with $w_i\geq 0$ can prove stability
of a constrained homogeneous random walk. It is easy to see that a
linear function $\Phi(q)=w^T\cdot q$ is a Lyapunov function if and
only if for some $\gamma>0$ and every nonempty set
$\Lambda\subset\{1,2,\ldots,d\}$ the following inequality holds

\be{LinearLyapunov} E[w^TQ(t+1)-w^TQ(t)|Q(t)\in{\cal
Z}_{\Lambda}]= \sum_{\Delta\in\{-1,0,1\}^d}(w^T\Delta)
p(\Lambda,\Delta)\leq -\gamma. \ee

The existence of a linear Lyapunov function is only sufficient but
not necessary for stability of the constrained random walk. It is
also useful sometimes to consider a geometric Lyapunov function,
defined as follows.

\begin{Defi}\label{GeometricLyapunov}  A function $\Phi_g:\cZ\rightarrow [1,+\infty)$ is defined
to be a geometric Lyapunov function with drift $0<\gamma_g<1$ and
exception set ${\cal B}\subset \cZ$ if $|{\cal B}|<\infty$ and for
every state $q\notin{\cal B}$ \be{drift_g}
{E[\Phi_g(Q(t+1)|Q(t)=q]\over\Phi_g(q)}=\sum_{q\in\cZ}{\Phi_g(q')\over\Phi_g(q)}p(q,q')\leq\gamma_g<1.
\ee
\end{Defi}

A geometric Lyapunov function is used, for example in Meyn and
Tweedie \cite{metwee_paper}, to prove exponentially fast mixing of
a Markov chain which admits a geometric Lyapunov function. The
precise statement of this result will be given below in Section
\ref{StationaryDecidable}. If the condition (\ref{LinearLyapunov})
is met for some $w$ and $\gamma$, then a function of the form
$\Phi_g(q)=\exp(\delta w^T\cdot q)$ is a geometric Lyapunov
function for a suitable value of $\delta>0$.

Throughout the paper we will assume all the states $q$ communicate
with the state $0$, that is $p^{(t)}(q,0)>0$ for some $t\geq 0$.
As a consequence, the random walk is irreducible. If it is in
addition positive recurrent, then it possesses a unique stationary
distribution $\pi:\cZ\rightarrow [0,1]$. Namely
$\sum_{q\in\cZ}\pi(q)=1$ and for any state $q$ \be{stationarity}
\sum_{q'\in\cZ}\pi(q')p(q',q)=\pi(q) \ee This stationary
distribution is defined completely by the set of transition
parameters $p(\Lambda,\Delta)$. Computing the stationary
probability distribution for these walks is the main focus of this
paper. It was established by the author in
\cite{gamarnik_decidability} that checking positive recurrence of
a constrained homogeneous random walk is an undecidable problem -
no algorithm can exist to achieve this task. However, if one is
lucky to construct a Lyapunov function, for example by checking
condition \ref{LinearLyapunov} for some nonnegative vector
$w\in\Re_+^d$, then the random walk is in fact positive recurrent.
One might be tempted to believe that in this case the analysis of
the random walk is simplified significantly. In Section
\ref{StationaryUndecidable} we show even if a linear Lyapunov
function exists, computing the stationary probability distribution
is still an undecidable problem. As in the case of stability
analysis, our main tool for establishing this undecidability
result is a counter machine and the halting problem defined in
Section \ref{CounterMachine}.

\section{Counter Machines, Halting Problem and Undecidability}\label{CounterMachine}
A counter machine (see \cite{blondel}, \cite{hopcroft}) is a
deterministic computing machine which is a simplified version of a
Turing machine -- a general description of an algorithm working on
a particular input. In his classical work on the Halting Problem,
Turing showed that certain decision problems simply cannot have a
corresponding solving algorithm, and thus are undecidable. For a
definition of a Turing machine and the Turing Halting Problem see
\cite{sipser}.  Ever since many quite natural problems, in
mathematics and computer science were found to be undecidable.
Some of the undecidability results in control theory were obtained
by reduction from a counter machine, see Blondel et al.
\cite{blondel}. For a survey of  decidability results in control
theory area see Blondel and Tsitsiklis  \cite{blondel_survey}.

A counter machine is described by 2 counters $R_1,R_2$ and a
finite collection of states $S$. Each counter contains some
nonnegative integer in its register. Depending on the current
state $s\in S$ and depending on whether the content of the
registers is positive or zero, the counter machine is updated as
follows: the current state $s$ is updated to a new state $s'\in S$
and one of the counters has its number in the register incremented
by one, decremented by one or no change in the counters occurs.

Formally, a  counter machine is a pair $(S,\Gamma)$.
$S=\{s_0,s_1,\ldots,s_{m-1}\}$ is a finite set of states and
$\Gamma$ is configuration update function $\Gamma:S\times
\{0,1\}^2\rightarrow S\times\{-2,-1,0,1,2\}$. A configuration of a
counter machine is an arbitrary triplet $(s,z_1,z_2)\in S\times
{\cal Z}_+^2$. A configuration $(s,z_1,z_2)$ is updated to a
configuration $(s',z_1',z_2')$ as follows. First a  binary vector
$b=(b_1,b_2)$ is computed were $b_i=1$ if $z_i>0$ and $b_i=0$ if
$z_i=0$, $i=1,2$. If $\Gamma(s,b)=(s',1)$, then the current state
is changed from $s$ to $s'$, the content of the first counter   is
incremented by one and the second counter does not change:
$z_1'=z_1+1, z_2'=z_2$. We will also write
$\Gamma:(s,z_1,z_2)\rightarrow (s',z_1+1,z_2)$ and
$\Gamma:s\rightarrow s',\Gamma:z_1\rightarrow
z_1+1,\Gamma:z_2\rightarrow z_2$. If  $\Gamma(s,b)=(s',-1)$, then
the current state becomes $s'$, $z_1'=z_1-1,z_2'=z_2$. Similarly,
if  $\Gamma(s,b)=(s',2)$ or  $\Gamma(s,b)=(s',-2)$, the new
configuration becomes $(s',z_1,z_2+1)$ or $(s',z_1,z_2-1)$,
respectively. If $\Gamma(s,b)=(s',0)$ then the state is updated to
$s'$, but the contents of the counters do not change. This
definition can be extended to the one which incorporates more than
two counters, but, in most cases, such an extension is not
necessary for our purposes.

Given an initial configuration $(s^0,z^0_1,z^0_2)$ the counter
machine uniquely determines subsequent configurations
$(s^1,z^1_1,z^1_2), (s^2,z^2_1,z^2_2),
\ldots,(s^t,z^t_1,z^t_2),\ldots \,\, .$ We fix a certain
configuration $(s^*,z_1^*,z_2^*)$ and call it a \emph{halting}
configuration. If this configuration is reached then the process
halts and no additional updates are executed. The following
theorem establishes the undecidability of the halting property.

\begin{theorem}\label{counter_undecidable}
Given a counter machine $(S,\Gamma)$, initial configuration
$(s^0,z_1^0,z_2^0)$ and the halting configuration
$(s^*,z_1^*,z_2^*)$, the problem of determining whether the
halting configuration is reached in finite time is undecidable. It
remains undecidable even if the initial and the halting
configurations are the same with both counters equal to zero:
$s^0=s^*,z_1^0=z_2^0=z_1^*=z_2^*=0$.
\end{theorem}

The first part of this theorem is a classical result and can be
founded in \cite{hooper}. The restricted case of $s^0=s^*,
z_i^0=z_i^*,i=1,2$ can be proven similarly by extending the set of
states and the set of transition rules. It is the restricted case
of the theorem which will be used in the current paper.

\section{Computing the stationary probability distribution. The undecidability result}\label{StationaryUndecidable}
Theorem \ref{counter_undecidable} was used in
\cite{gamarnik_decidability} to prove that the stability of a
constrained random walk in $\cZ$ is undecidable. Naturally, the
problem of stability comes before the problem of computing the
stationary distribution of a stable random walk. As we mentioned
in Section \ref{definitions}, stability can be checked sometimes
by constructing a Lyapunov function. In this section we prove our
main result: even if such a Lyapunov function, witnessing
stability, is available and is provided as a part of the data
parameters, computing stationary distribution is an undecidable
problem.

We now give an informal outline of the proof. The proof  uses a
reduction from a halting problem for a counter machine. We embed a
counter machine with initial and halting configuration $(s^*,0,0)$
into a deterministic walk in $\cZ$. The state space and the
transition rules of this walk are then extended in some way that
incorporates an independent Bernoulli process with some fixed
parameter $p$. We then show that
\begin{itemize}
\item If the original counter machine never returns to the
initial configuration $(s^*,0,0)$, then the constructed random
walk, when started from the origin, returns into the origin in
$2t+2$ steps with probability $(1-p)p^t$ for $t=0,1,2,\ldots\,\,.$
In particular, the expected return time to the origin is
$2/(1-p)$.
\item If the original counter machine returns to the initial
configuration in $T$ steps, then the modified random walk returns
into the origin in $2t+2$ steps with probability $(1-p)p^t$ for
$t\leq T-1$ and in $2T+2$ steps with the remaining probability
$1-\sum_{t\leq T-1}(1-p)p^t=p^T$. In particular, the expected
return time to the origin is $(2-2p^{T+1})/(1-p)$.
\end{itemize}

The stationary probability distribution of any state is exactly
the expected return time to this state. Therefore, the stationary
probability of the origin is $(1-p)/2$ if the counter
machine halts and is strictly greater, if the counter machine does
not halt. Since the value $p$ is our control parameter, and since
checking whether the counter machine halts is an undecidable
problem, then computing stationary probability is undecidable as
well. We now state and prove rigorously this result. As before,
let $\pi$ denote the unique stationary distribution of an
irreducible positive recurrent random walk. Let also $0$ denote
the origin of the nonnegative lattice $\cZ$.

\begin{theorem}\label{MainResultRW}
Given an irreducible constrained random walk with transition
probabilities $p(\Lambda,\Delta)$, given a linear vector $w\in
\Re_+^d$ satisfying (\ref{LinearLyapunov}) and given a rational
value $0\leq r\leq 1$, the problem of checking whether $\pi(0)\leq
r$ is undecidable -- no algorithm exists which achieves this task.
\end{theorem}

\remarks

1. The stationary distribution can in principle take non-rational
values. In order to put the problem into a framework suitable for
algorithmic analysis we modified the question into the one of
checking whether $\pi(\cdot)\leq r$ for rational values $r$. This
is a standard method in the theory of Turing decidable numbers,
see \cite{RealComplexity}.

2. A simple example where computing the stationary probability
distribution is a decidable problem is  Jackson networks,
\cite{kleinrock}. For such a network with $d$ stations the
stationary probability of the state $m=(m_1,\ldots,m_d)$ is given
by $\prod_{j=1}^d(1-\rho_j)\rho_j^{m_j}$, where $\rho_j$ is the
traffic intensity in station $j$. Specifically, the stationary
probability of the state $0$ is $\prod_{j=1}^d(1-\rho_j)$. Given
any rational value $0\leq r\leq 1$, it is a trivial computation to
check whether this product is at least $r$.

\vvert {\it Proof of Theorem \ref{MainResultRW}:} we start with a
construction used  in \cite{gamarnik_decidability}. Namely, we
embed a given counter machine with states $s_0,s_1,\ldots,s_{m-1}$
into a deterministic walk in $\cZm$ as follows. Without the loss
of generality, assume that $s^*=s_0$. Let configuration
$(s_i,z_1,z_2), 1\leq i\leq m-1 $ correspond to the state
$q=(e_i,z_1,z_2)\in\cZm$, where $e_i$ is unit vector with $1$ in
$i$-th coordinate and zero everywhere else. Also, let
configurations $(s_0,z_1,z_2)$ correspond to $(0,z_1,z_2)$, with
zeros in first $m-1$ coordinates. Specifically, the initial and
halting configuration $(s_0,0,0)$ corresponds to the origin $0$.
We now describe the set of transition vectors
$\Delta=\Delta(\Lambda)$. We describe it first for subsets
$\Lambda\subset\{1,2,\ldots,m+1\}$ which correspond to an encoding
of some configuration of a counter machine. Specifically,
$\Lambda\cap\{1,2,\ldots,m-1\}=\emptyset$ (corresponding to
configurations with state $s_0$) or
$\Lambda\cap\{1,2,\ldots,m-1\}=\{i\}$ for some $1\leq i\leq m-1$,
corresponding to configurations with state $s_i$. Fix any
configuration $(s_i,z_1,z_2)$. Suppose the corresponding update
rule is $\Gamma((s_i,z_1,z_2))=(s_j,+1)$ for some $1\leq j\leq
m-1$. That is, the state is changed into $s_j$, the first counter
is incremented by $1$ and second counter remains unchanged. We
make the corresponding transition vector to be
$\Delta=\Delta(\Lambda)$, where the $i$-th coordinate of $\Delta$
is $-1$, the $j$-th coordinate is $+1$, the $m$-th coordinate is
$+1$ and all the other coordinates are zeros. It is easy to see
that if at time $t$, the state $Q(t)$ corresponds to some
configuration $(s_i,z_1,z_2)$, that is $Q(t)=(e_i,z_1,z_2)$, then
$Q(t+1)=Q(t)+\Delta$ corresponds to the configuration
$(s_j,z_1+1,z_2)$ obtained by applying rule $\Gamma$. We construct
transition vectors similarly for other cases of configuration
updates. In particular, if the state $s_i$ is changed to state
$s_0$, then the corresponding $\Delta$ has $-1$ in the $i$-th
coordinate and zeros in all the coordinates $1\leq j\leq m-1,j\neq
i$. As we will see later, if $Q(t)$ corresponds to some
configuration of a counter machine at time $t$, then it does so
for all the later time $t'\geq t$. Now if $Q(t)$ belongs to some
face ${\cal Z}_{\Lambda}$ which does not correspond to some
configuration, then we simply set $\Delta(\Lambda)=-e_i$ where $i$
is the smallest coordinate which belongs to $\Lambda$. Then at
some later time $t'>t$ the state $Q(t')$ will correspond to some
configuration.

Construction above is exactly the one used in
\cite{gamarnik_decidability} to analyze stability. We now modify
the construction by adding two additional coordinates. Our new
state at time $t$ is thus denoted by $\bar
Q(t)=(Q(t),q_1(t),q_2(t))\in \cZmm$. Also a parameter $0<p<1$ is
fixed. The transition rules are modified as follows.
\begin{enumerate}
\item When $q_2(t)=1$, the first part $Q(t)$ of the state is
updated exactly as above. Also, if $||Q(t)||>0$, in other words,
$Q(t)$ does not represent the halting configuration $(s_0,0,0)$,
then the value of $q_2(t)$ stays $1$ with probability $p$ and
switches to $0$ with probability $1-p$. If, on the other hand
$||Q(t)||=0$ then we set $q_2(t+1)=0$ with probability $1$.
Finally, the value of $q_1(t)\in\{-2,-1,0,1,2\}$ is selected in
such a way that $||(Q(t+1),q_1(t+1))||=||(Q(t),q_1(t))||+1$, where
$||(Q(t),q_1(t))||=\sum_{i=1}^{m+1}Q_i(t)+q_1(t)$. It is  easy to
see that such a value of $q_1(t)$ always exists. For example if
$Q(t)$ encodes $(s_i,z_1,z_2),i\neq 0$ and the configuration is
changed into $(s_j,z_1,z_2-1),j\neq 0$, then we put $q_1(t)=2$.

\remark We stipulated before that the transition vectors $\Delta$
must belong to $\{-1,0,1\}^{m+4}$ for our constrained random walk,
whereas above the value of $q_1(t)$ can change by  $-2$ and $2$.
It is easy to satisfy this constraint for $q_1(t)$ by splitting it
into two coordinates $q_1(t),q_1'(t)$ and making
$q_1(t)=q_1'(t)=1$ in case $q_1(t)$ was assigned $2$ before, and
$q_1(t)=q_1'(t)=-1$ in case $q_1(t)$ was assigned $-2$. We keep
only one $q_1(t)$ for simplicity, allowing it to take values
$-2,2$.

\item When $q_2(t)=0$, we set $\Delta_k=-1,\Delta_i=0,i\neq k,1\leq i\leq m+3$,
where $k$ is the smallest coordinate such that $\bar Q_k(t)>0$. In
particular, $q_2(t)$ stays equal to $0$. If $Q(t)=q_1(t)=0$, (in
particular $Q(t)$ encodes the initial-terminal configuration
$(s_0,0,0)$) then $Q(t)$ and $q_1(t)$ are updated as in the case
$q_2(t)=1$ above. Also $q_2(t)$ in this case is switched to $1$
with probability $p$ and stays $0$ with probability $1-p$.
\end{enumerate}

Note, that the only stochastic part in our random walk is  the
last component $q_2(t)$.

\begin{prop}\label{RecurrenceTime}
The constructed random walk $\bar Q(t)$ is irreducible and
positive recurrent with the unique stationary distribution $\pi$.
Moreover,
\begin{enumerate}
\item If the counter machine with the initial configuration
$(s_0,0,0)$ does not halt, then the random walk $\bar Q(t)$ with
the initial state $\bar Q(0)=0$ returns to the origin in $2t+2$
steps with probability $(1-p)p^t$, for $t=0,1,2,\ldots\,\,.$ As a
result, the expected recurrence time of the state $0$ is
$1/\pi(0)=2/(1-p)$.
\item If the counter machine with the initial configuration
$(s_0,0,0)$ halts in $T\geq 1$ steps, then the random walk $\bar
Q(t)$ with the initial state $\bar Q(0)=0$ returns to the origin
in $2+2t$ steps with probability $(1-p)p^t$ for $t<T$, and in
$2+2T$ steps with the remaining probability $p^T$. As a result,
the expected recurrence time of the state $0$ is
$1/\pi(0)=(2-2p^{T+1})/(1-p)$.
\item For any $C\geq 2/(1-p)$ the function
$\sum_{i=1}^{m+1}Q_i+q_1(t)+Cq_2(t)$ is a linear Lyapunov function
with drift $\,-\gamma=-1$ and an exception set ${\cal B}=\{0 \}$.
\end{enumerate}
\end{prop}

We first show that the proposition above implies the theorem.
Suppose, we had an algorithm ${\cal A}$ which given an irreducible
constrained random walk $Q(t)$, with a linear Lyapunov function
$w^TQ(t)$ and given a rational value $0\leq r\leq 1$ could
determine whether the unique stationary distribution $\pi$
satisfies $\pi(0)\leq r$. We take a counter machine and construct
a random walk $\bar Q(t)$ as described above. Proposition
\ref{RecurrenceTime} implies that this walk is a valid input for
the algorithm ${\cal A}$. We use ${\cal A}$ to determine whether
$\pi(0)\leq r\equiv (1-p)/2$. From Proposition
\ref{RecurrenceTime}, this is the case if and only if the
underlying counter machine does not halt. In this fashion, we
obtain an algorithm for checking halting property for counter
machines. This is a contradiction to Theorem
\ref{counter_undecidable}. \qed

\vvert {\it Proof of Proposition \ref{RecurrenceTime}:} Suppose
the underlying counter machine does not halt. Let us trace the
dynamics of our random walk $\bar Q(t)$ starting from $\bar
Q(0)=0$. Initially, by applying rule 2,  it moves into some state
$(Q(1),0,1)$ with probability $p$ or state $(Q(1),0,0)$ with
probability $1-p$. An independent Bernoulli process for $q_2(t)$
with parameter $p$ is continued in the first case. Suppose this
process succeeds exactly $t\geq 0$ times (including the transition
from initial state $0$), which occurs with probability $(1-p)p^t$.
Then, applying rule 1, at times $t$ and $t+1$ we have states
$(Q(t),q_1(t),1),(Q(t+1),q_1(t+1),0)$ with
$||Q(t)+q_1(t)||=t,||Q(t+1)+q_1(t+1)||=t+1$. At this moment rule 2
becomes applicable. Since at each step the norm $||\bar Q(t)||$
decreases exactly by one, the origin is reached at time
$(t+1)+(t+1)$. We conclude that the return time is $2+2t$ with
probability $(1-p)p^t, t=0,1,2,\ldots\,\,$. The expected return
time is then $2/(1-p)$ and the stationary probability of the state
$0$ is $(1-p)/2$.

Suppose, now, the underlying counter machine reaches the terminal
state $(s_0,0,0)$ in exactly $T\geq 1$ steps. Suppose also the
Bernoulli process for $q_2(t)$ succeeds exactly $t\geq 0$ times
for $t<T$. Then, exactly as above, the origin is reached in $2+2t$
steps and this occurs with probability $(1-p)p^t$. If, however, by
the time $T$ the Bernoulli process does not fail, which occurs
with probability $p^T$, then the state $\bar Q(T)=
(Q(T),q_1(t),1)=(0_{m+1},q_1(t),1)$ is reached at time $T$, where
$0_{m+1}$ denotes a $m+1$-dimensional zero vector. By the choice
of $q_1(t)$ in rule 1, $q_1(T)=T$. At time $T+1$, by rule 1, we
have a state $\bar Q(T+1)=(0_{m+2},T+1,0)$ and rule 2 applies. At
time $(T+1)+(T+1)$ the origin is reached. We conclude that the
random walk returns to the origin in $2+2T$ steps. Combining, the
expected return time to the origin is then $(2-2p^{T+1})/(1-p)$,
if the counter machine halts in $T$ steps, and the stationary
probability of the state $0$ is $(1-p)/(2-2p^{T+1})$.

To complete the proof of the proposition, we analyze the expected
change of the function $\sum_{i=1}^{m+1}Q_i+q_1(t)+Cq_2(t)$. When
$q_2(t)=0$ and $\bar Q(t)\neq 0$, the sum decreases
deterministically by $1$. When $q_2(t)=1$, the value of
$\sum_{i=1}^{m+1}Q_i+q_1(t)$ increases deterministically by $1$,
and the value of $Cq_2(t)$ stays the same with probability $p$ or
decreases by $C$ with probability $1-p$. Therefore, the expected
change of the sum is $1-C(1-p)$. When $C\geq 2/(1-p)$, the
expected change is at most $-1$. \qed

\vvert

An important implication of Theorem \ref{MainResultRW} is that it
is impossible to express the stationary distribution of a positive
recurrent random walk $Q(t)$ as a function of the parameters
$p(\Lambda,\Delta)$ via some computable function $f(\cdot)$. For
example, the stationary distribution cannot be expressed as roots
of some polynomial equations with rational coefficient, as
inequalities $x\leq r$ can be checked for any root $x$ of such a
polynomial and any rational value $r$. This is a startling
contrast to a simple expression $\prod (1-\rho_j)\rho_j^{n_j}$
corresponding to a stationary distribution of a Jackson network.

\section{Large Deviation Rates. The undecidability
result}\label{LDUndecidable}

In this section we discuss the question of computing large
deviation rates for our model. Specifically, we focus on computing
large deviation rates for the stationary distribution $\pi$ of our
random walk $Q(t)$ in $\cZ$. Let $\Re$ and $\Re_+$ denote the set
of real values and the set of nonnegative real values,
respectively. For any $x\in \Re$ let $\lfloor x\rfloor$ denote
largest integer not bigger than $x$, and for any $x\in\Re^d$ let
$\lfloor x\rfloor=(\lfloor x_1\rfloor,\ldots,\lfloor x_d\rfloor)$.
We say that a function $L:\Re_+^d\rightarrow \Re_+\cup\{\infty\}$
is a large deviation rate function for a given irreducible
positive recurrent random walk $Q(t)$ in $\cZ$ if for any vector
$v\in\Re_+^d$, the stationary distribution $\pi$ satisfies

\be{ldrates} \lim_{n\rightarrow\infty}{\log(\pi(\lfloor
vn\rfloor))\over n}=L(v), \ee

In other words, the stationary probability of being in state
$\lfloor vn\rfloor$ is asymptotically $\exp(-L(v)n)$ for large
$n$. For results on large deviations for specific types of
constrained random walks in $\cZ$ see
\cite{IgnatyukMalyshevScherbakov}. There are numerous works on
large deviation in the context of queueing systems, see Shwartz
and Weiss \cite{large_deviations} for a survey. Specifically,
Kurkova and Suhov \cite{suhov} study large deviation rates for a
two dimensional random walk corresponding to join-the-shortest
queue. The analysis is quite intricate and uses complex-analytic
techniques developed by Malyshev \cite{malyshev},
\cite{MalyshevAnalyticMethods}, \cite{MalyshevAsymptotic} back in
70's. To the best of our knowledge, the existence of the large
deviations limits (\ref{ldrates}) is not fully proved for general
constrained homogeneous random walks $Q(t)$ in $\cZ$. One can
instead consider limits

\be{ldrates+-}
L_{-}(v)\equiv\liminf_{n\rightarrow\infty}{\log(\pi(\lfloor
vn\rfloor))\over n},\qquad
L_{+}(v)\equiv\limsup_{n\rightarrow\infty}{\log(\pi(\lfloor
vn\rfloor))\over n}. \ee

The goal of the present section is to prove that computing the
large deviation rate function $L(v)$ is an undecidable problem,
even if the walk is known to be a priori positive recurrent via,
for example, existence of a linear Lyapunov function, and even if
the large deviation limit function $L(v)$ is known to exist. The
following is the main result of this section.

\begin{theorem}\label{ldr}
Given an irreducible constrained random walk, given a linear
vector $w\in \Re_+^d$ satisfying (\ref{LinearLyapunov}), given a
rational value $0\leq r\leq 1$ and a vector $v\in Z_+^d$, the
problems of determining whether $L_{-}(v)\leq r,L_{+}(v)\leq r$
are undecidable.
\end{theorem}

\remark As we will see below, the large deviations limit function
$L(v)=L_{-}(v)=L_{+}(v)$ exists for the subclass of random walks
we consider. As before, the reason for including a linear Lyapunov
function into the condition of the theorem is to provide a simple
way of insuring that the walk is positive recurrent.

\vvert \proof The proof is again based on reduction from a halting
problem for a counter machine. Given a counter machine with $m$
states consider the extended $m+3$-dimensional random walk $\bar
Q(t)$ constructed in the proof of Theorem \ref{MainResultRW}. We
extend it even further by adding an additional coordinate $(\bar
Q(t), q_3(t))$. Recall, that the rules for updating $q_1(t)$ were
such that $||Q(t)||+q_1(t)=t$ as long as $q_2(t')=1$ for $1\leq
t'\leq t$. Construct the rules for updating $q_3(t)$ as follows.
If $q_2(t)=1$, then $q_3(t+1)=q_3(t)+1$. Also, if
$Q(t)=q_1(t)=q_2(t)=0$ then again $q_3(t+1)=q_3(t)+1$. In other
words, as long as the random walk starts from the origin and as
long as $q_2(t)$ remains equal to $1$, $q_3(t)=||Q(t)||+q_1(t)=t$.
Once $q_2(t)$ becomes zero, the value of $q_3(t)$ stays the same
as long as $\bar Q(t)\neq 0$ and decreases by one when $\bar Q(t)$
becomes zero, and continues decreasing until it itself becomes
zero.

Let $v=(0,\ldots,0,1)$ be an $m+4$ dimensional vector with the
last coordinate equal to unity and all other coordinates equal to
zero. We now analyze the large deviation rate $L(v)$ for this
vector with respect to the unique stationary distribution $\pi$.
Specifically, we show that the value of $L(v)$ depends on whether
the counter machine halts. Indeed, if the counter machine halts in
$T$ steps, then the value of $||Q(t)||+q_1(t)+q_2(t)+q_3(t)$ is
bounded by $2T+1$ and, as a result, the stationary probability of
the state $nv$, for large $n$ becomes zero. That is
$L(v)=+\infty$. Now we show that if the counter machine does not
halt, then $L(v)=\log p$. We compute $\pi(nv)$ by computing the
expected return time to state $nv$, when the random walk is in
this state at time $0$. Thus, we have $q_3(0)=n,
Q(0)=q_1(0)=q_2(0)=0$. By the update rules of $q_3$, it decreases
by one at each time step and at time $t=n$ it becomes zero. All
the other components remain equal to zero. Beginning from this
time $t=n$, the random walks keeps returning to the origin $0$
after some random time intervals. It is easy to see that the
probability that the state $nv$ is visited in between any given
two visits to the origin is exactly $p^n$ -- the probability that
the Bernoulli process survives at least $n$ steps. Let
$R_1,R_2,\ldots\,\,$ denote the random time intervals between
successive visits to the origin. For a fixed $m\geq 1$ the
probability that $R_m$ is the first interval during which $vn$ is
visited is $p^n[1-p^n]^{m-1}, m=1,2,\ldots$, and the expected
number of intervals $R_m$ before state $nv$ is visited for the
first time is $(p^n)^{-1}$. Let $I_n$ denote the indicator
function for the event "state $vn$ is visited between visits to
the origin". In particular, $\pr\{I_n\}=p^n$. We now compute
$E[R_m|I_n]$ and $E[R_m|\bar I_n]$. Note, that the time $n$ which
takes to get from $nv$ to the origin plus the expected time  it
takes to get from the origin to $nv$, conditioned on $I_n$ is
exactly $E[R_m|I_n]$. We then obtain that the expected recurrence
time of the state $nv$ is $(p^n)^{-1}E[R_m|\bar I_n]+E[R_m|I_n]$.

To compute $E[R_m|\cdot]$, recall from Proposition
\ref{RecurrenceTime} that

\be{Rm} E[R_m]=E[R_m|I_n]\pr\{I_n\}+E[R_m|\bar I_n]\pr\{\bar
I_n\}=2/(1-p).
\ee

If the Bernoulli process survives $t\geq n$ steps then $R_m=3t+1$
steps. Then

\be{In} E[R_m|I_n]={E[R_mI_n]\over \pr\{I_n\}}={\sum_{t\geq
n}(1-p)p^t(3t+1)\over p^n}={(4-p)p^n\over p^n}=4-p
\ee

Then we obtain from (\ref{Rm})
\[
E[R_m|\bar I_n]={2/(1-p)-(4-p)p^n\over 1-p^n}
\]
We conclude that the expected return time to the state $nv$ is
\[
{1\over \pi(nv)}=(p^n)^{-1}{2/(1-p)-(4-p)p^n\over 1-p^n}+(4-p),
\]
and
\[
\lim_{n\rightarrow\infty}{\log(\pi(vn))\over n}=\log p
\]
as we claimed. We see that the value of $L(v)$ depends on whether
the underlying counter machine halts or not. Specifically, by
taking any rational value $r>\log p$, we conclude that the problem
of checking whether $L(v)\leq r$ is undecidable, by appealing
again to Theorem \ref{counter_undecidable}. \qed

\remarks

1. Note that we cannot determine the value of $L(v)$ even
approximately, as we cannot distinguish between the cases
$L(v)<+\infty$ and $L(v)=+\infty$. Contrast this with the results
of Section \ref{StationaryDecidable}.

2. We would not need extra coordinate $q_3(t)$ if we were
interested in large deviation rate \\
$\lim_{n\rightarrow\infty}\pi(||\bar Q(t)||)/n$ of the stationary
distribution of the norm of the state. The analysis would be
identical to the one above.

\section{Application to queueing systems}\label{queues}
The results of the previous sections have implications to a
certain type of queueing systems. A queueing system consisting of
a single station processor and operating under a certain class of
{\it generalized priority} policies was introduced in
\cite{gamarnik_decidability}. It was shown that, similar to
constrained random walks, determining stability for these queueing
systems is an undecidable problem. In this section we consider the
same class of system and show that computing stationary
probabilities and large deviation rates are undecidable problems
as well.

We start with the description of the system. Consider a single
station queueing system ${\cal Q}$ consisting of a single server
and $I$ types of parts arriving externally. The parts
corresponding to type $i=1,2,\ldots,I$ visit the station $J_i$
times. On each visit each part must receive a service before
proceeding to the next visit. Only one part among all the types
can receive service at a time. While waiting for service for the
$j$-th time, the type $i$ part is stored in buffer $B_{ij}$. We
denote by $n$ the total number of buffers $n=\sum_{i=1}^IJ_i$. The
service time for each part in each visit is assumed to be equal to
unity. Each part can arrive into the system only in times which
are  multiples of some fixed integer value $M$. Specifically,
certain values $0\leq p_i\leq 1$ are fixed for each type $i$. For
each type $i$ and each $m=0,1,2,\ldots,$ exactly one part arrives
at time $mM$ with probability $p_i$ and no part arrives with
probability $1-p_i$, independently for all $m$ and all other
types. In particular, interarrival times are geometrically
distributed with expected interarrival time equal to
$1/\lambda_i=M/p_i$, where, correspondingly, $\lambda_i$ is the
arrival rate for type $i$.

A scheduling policy $u$ is defined to be a {\it generalized
priority} policy if it operates in the following manner. A
function $u:\{0,1\}^n\rightarrow \{0,1,2,\ldots,n\}$ is fixed. At
each time $t=0,1,2,\ldots $ the scheduler looks at the system and
computes the binary vector $b=(b_1,b_2,\ldots,b_n)\in\{0,1\}^n$,
where $b_i=1$ if there are parts in the $i$-th buffer and $b_i=0$,
otherwise. Then the value $k=u(b), 0\leq k\leq n$ is computed. If
$k>0$ then the station processes a part in the $k$-th buffer. If
$k=0$ the server idles. The map $u$ is assumed to satisfy the
natural consistency condition: $u(b)=k>0$ only when $b_k=1$. That
is, processing can be done in buffer $k$ only when there are jobs
in buffer $k$. Note that the generalized priority scheduling
policy is defined in finitely many terms and is completely state
dependent - the scheduling decision at time $t$ does not depend on
the state of the queueing system at times $t'<t$. A usual priority
policy corresponds to the case when there is some permutation
$\theta$ of the buffers $\{1,2,\ldots,n\}$ and $u(b)=k$ if and
only if $b_k=1$ and $b_i=0$ for all $i$ such that
$\theta(i)<\theta(k)$. In words, priority scheduling policy
processes parts in buffers with lowest value (highest priority)
$\theta$, which still has parts. Once we specify the queueing
system ${\cal Q}$ and some generalized priority policy $u$ we have
specified some discrete time discrete space stochastic process.
This process considered in times $t=mM,m=0,1,2,\ldots\,\,$ is in
fact a Markov chain.

Given a generalized priority policy $u$,  a pair $({\cal Q},u)$ is
defined to be stable if there exists a finite number $C>0$ such
that the total number of parts in the queueing system ${\cal Q}$
at time $t$ does not exceed $C$ for infinitely many $t$ with
probability $1$. In other words, the underlying Markov chain is
positive recurrent. In this case there exists at least one
stationary probability distribution. It is known that the
necessary condition for stability is the following load condition
\begin{equation}\label{load}
\rho\equiv \sum_{i=1}^I\sum_{j=1}^{J_i}\lambda_i<1.
\end{equation}
This condition is  also sufficient for stability if  the policy is
work conserving, which does not apply here, since we allow idling
$u(b)=0$. We assume that the load condition above holds. We define
a Lyapunov function and large deviations rates $L(v)$ for this
queueing system in the same way we did for constrained homogeneous
random walks in Section \ref{LDUndecidable}. As for constrained
random walks, we now show that computing  stationary probability
distributions and computing large deviations rates for queueing
systems operating under  generalized priority policies is not
possible. As for constrained random walks, we show that these
problems are impossible to solve even if the underlying Markov
chain is known to be irreducible and a linear Lyapunov function is
available. Let $\pi$ denote the unique stationary distribution of
a given irreducible positive recurrent queueing system $({\cal
Q},u)$. Let also $0$ denote the state of the system with all
buffers empty.

\begin{theorem}\label{QueueUndecidable}
Given a queueing system ${\cal Q}$ operating under some
generalized priority policy $u$, given a linear Lyapunov function
$\Phi$ and given a rational value $0\leq r\leq 1$, the problem of
determining whether $\pi(0)\leq r$ is undecidable. Likewise, given
a vector $v$, the problem of determining whether $L(v)\leq r$ is
undecidable.
\end{theorem}

\vvert \proof A reduction from a counter machine to a queueing
system operating under some generalized priority policy was
constructed in \cite{gamarnik_decidability}. This reduction had
the following features. Given a counter machine with $m$ states,
the corresponding queueing system had $24$ buffers and $I=3m+7$
arrival streams. There is a one-to-one correspondence between the
configurations of the counter machine and states of the queueing
system.  In particular, if a counter machine has configuration
$(s_i,z_1,z_2)$ at time $t$, then the queueing system at time
$(3m+26)t$ is in a state which corresponds to this configuration
in some well-defined way. We omit the details of this reduction
and instead refer the reader to \cite{gamarnik_decidability}. We
now modify the reduction to incorporate the extended random walk
$\bar Q(t)=(Q(t),q_1(t),q_2(t))$ that was constructed in Section
\ref{StationaryUndecidable}. Recall, that the part $Q(t)$ of this
walk represented exactly $m$ states and the two counters of the
underlying counter machine. We add two additional streams of
arrivals which correspond to coordinates $q_1(t)$ and $q_2(t)$. We
also construct additional buffers for $q_1$ and $q_2$ exactly in
the way we did in \cite{gamarnik_decidability} for counters
$z_1,z_2$. The interarrival  times for all the arrival streams,
except for the stream corresponding to $q_2$, are  deterministic
and equal to some integer $M$ which is selected to be bigger than
the number of buffers. For the stream corresponding to $q_2$, at
most one part arrives at times $Mt,t=0,1,2,\ldots$ independently
for all $t$, and the probability that a part does arrive at time
$t$ is equal to $p$, where $p$ is the parameter selected in
construction of the random walk $\bar Q(t)$. Thus, $p_i=p$ for the
arrival streams corresponding to $q_2$ and $p_i=1$ for all the
other arrival streams. Finally, we modify the rules of the
generalized priority policy to incorporate the rules by which the
values of $q_1(t),q_2(t)$ are updated. This can be done in a way
similar to the rules corresponding to $z_1,z_2$ in
\cite{gamarnik_decidability}. We thus obtain a system which mimics
the dynamics of $\bar Q(t)$ at times $Mt, t=0,1,2,\ldots\,\,$. A
linear Lyapunov function can be constructed again, provided that
the parameter $p$ is sufficiently small. Arguing as in the proof
of Theorem \ref{MainResultRW}, we show that the problem of
checking whether $\pi(0)\leq r$ is undecidable. Similarly, we show
that the problem of checking whether $L(v)\leq r$ is undecidable,
where $v$ is the unit vector with one in the coordinate
corresponding to $q_2(t)$ and zero in all the other coordinates.
For the latter case of computing large deviations rates, we add an
additional arrival stream and buffers to represent the part $q_3$.
\qed

\section{Computing stationary probabilities approximately using a
Lyapunov function}\label{StationaryDecidable}

In this section we show that, despite the results of Section
\ref{StationaryUndecidable}, computing the stationary probability
is possible, if we are willing to tolerate some two-sided error
and a computable geometric Lyapunov function $\Phi_g$ exists.  Our
result is  a simple consequence of the following result
established by Meyn and Tweedie \cite{metwee_paper}, which shows
that infinite state Markov chain  mixes exponentially fast when a
geometric Lyapunov function can be constructed. The following is
Theorem 2.3 proven in \cite{metwee_paper}.

\begin{theorem}\label{mixing}
Given an irreducible Markov chain $Q(t)$, suppose $\Phi_g$ is a
geometric Lyapunov function with a geometric drift $\gamma_g<1$
and the exception set ${\cal B}$. Suppose also that $\pi$ is the
unique stationary distribution. Then, there exist constants
$R>0,0<\rho<1$ such that for any state $x\in{\cal X}$ and any
function $\phi:{\cal X}\rightarrow \Re$ satisfying $\phi(x)\leq
\Phi(x),\forall\,\,x\in{\cal X}$, the following bound holds

\be{mixing_bound} \Big |\sum_{y\in{\cal
X}}\phi(y)\Big(\pr\{Q(t)=y|Q(0)=x\}-\pi(y)\Big)\Big |\leq
\Phi_g(x)R\rho^t. \ee

The constants $R,\rho$ are computable functions which depend on
$\gamma_g,\max_{x\in {\cal B}}\Phi_g(x)$ and

\be{pmin} \nu^{\Phi}_g=\max_{x,x'\in\cZ}\{{\Phi_g(x')\over
\Phi_g(x)}:p(x,x')>0\},\qquad p^{\cal B}_{\min}\equiv \min_{x,y\in
{\cal B}}p(x,y). \ee
\end{theorem}

Exact formulas for computing $R,\rho$ are provided in
\cite{metwee_paper}. They are quite lengthy and we do not repeat
them here. These formulas give meaningful bounds only in case
$0<\gamma_g<1;\nu^{\Phi}_g<\infty; p^{\cal B}_{\min}>0$.

Given a fixed state $x_0\in\cZ$, consider the function
$\phi(x_0)=1/\Phi_g(x_0), \phi(x)=0,x\neq x_0$. This function
satisfies the conditions of the theorem and one obtains a
computable bound on the difference
$|\pr\{Q(t)=x_0|Q(0)=x\}-\pi(x_0)|$, which decreases exponentially
fast with $t$. This bound can be used for computing stationary
probability distribution $\pi$.

\begin{theorem}\label{StationaryDecidableTheorem}
Given a constrained random walk $Q(t)$ in $\cZ$, given a state
$x_0\in\cZ$ and an arbitrary value $\epsilon>0$, under the
conditions of Theorem \ref{mixing}, there exists an computable
value $\hat x$ which satisfies $\pi(x)\in [\hat x-\epsilon,\hat
x+\epsilon]$. In other words, the stationary probability of the
state $x_0$ can be computed approximately with an arbitrary degree
of accuracy.
\end{theorem}

\vvert
\proof The proof is a simple consequence of Theorem
\ref{mixing}. We fix an arbitrary initial state $Q(0)$, say
$Q(0)=0$. Compute the values $R,\rho,1/\Phi_g(x_0)$. Select $t$
large enough, so that $\Phi_g(Q(0))R\rho^t<\epsilon$. Compute the
transient probability $Q(t)=x_0$ conditioned on $Q(0)=0$. This can
be done by direct calculation since $t$ is finite and from any
state there are only finitely many neighboring states that can be
entered with positive probability. The value
$\pr\{Q(t)=x_0|Q(0)=0\}$ can be taken as $\hat x$, using
inequality (\ref{mixing_bound}) and by the choice of $t$. \qed

As we mentioned above, a similar result cannot be established for
large deviations rates $L(v)$, since the value of $L(v)$ changes
between $\log p$ and $+\infty$ depending on whether the underlying
counter machine halts or not. Therefore, computing the value of
$L(v)$ even approximately still is an undecidable problem.

\section{Conclusions}\label{conclusions}
We considered in this paper the problems of computing stationary
probability distributions and large deviations rates for
constrained homogeneous random walks in $\cZ$. Both problems were
shown to be undecidable -- no algorithmic procedure for solving
these problems can exist. An implication of these results is that
no useful formulas for computing these  quantities, for example
along the lines of formulas for product form networks, can exist.
For the problems of computing stationary probabilities, we showed
that an approximate computation is possible with arbitrary degree
of accuracy if a suitable geometric Lyapunov function can be
constructed. Yet the problem of computing large deviation rates
remains to be undecidable even in approximation sense as even
checking whether a large deviation rate along a given vector is
finite or not, is an undecidable problems.

We conjecture that these problems remain to be undecidable in more
restrictive and interesting class of Markov chains corresponding
to multiclass queueing networks operating under more conventional
scheduling policies like First-In-First-Out or priority polices.


\providecommand{\bysame}{\leavevmode\hbox to3em{\hrulefill}\thinspace}
\providecommand{\MR}{\relax\ifhmode\unskip\space\fi MR }
\providecommand{\MRhref}[2]{%
  \href{http://www.ams.org/mathscinet-getitem?mr=#1}{#2}
}
\providecommand{\href}[2]{#2}

\end{document}